\documentstyle[12pt]{article}
\newcommand{\beq}{\begin{equation}}
\newcommand{\eeq}{\end{equation}}
\newcommand{\barr}{\begin{array}}
\newcommand{\earr}{\end{array}}
\newtheorem{theorem}{Theorem}
\newtheorem{lemma}{Lemma}

\newtheorem{remark}{Remark}
\newcounter{one}
\input amssym.def
\begin{document}
\begin{center}

{\Huge Closed curves in ${\bf R}^3$: a characterization in terms of
  curvature and torsion, the Hasimoto map and periodic solutions of
  the Filament Equation.}

\vspace{1cm}

{\bf P.G.Grinevich}\footnote{This work was completed during the
author's visit to the Freie Universit\"at, Berlin, Germany, which 
was supported by the Humboldt-Foundation. Part of this work was
fulfilled during the author's visit to the Maryland University. 
This work was also partially supported by the Russian Foundation for 
Basic Research, grant No 95-01-755. }

\medskip

{\bf M.U.Schmidt}\footnote{Supported by the DFG, SFB 288
``Differentialgeometrie und  Quiantenphysik''.}

\end{center}

\vspace{0.5cm}

${\vphantom k}^{\rm1}$ L.D.Landau Institute for Theoretical Physics, 
Kosygina 2, Moscow, 177940, Russia.\\ e-mail: pgg@landau.ac.ru

${\vphantom k}^{\rm2}$ Freie Universit\"at, Berlin, Arnimallee 14, D-14195,
Berlin, Germany.\\e-mail: mschmidt@physik.fu-berlin.de

\begin{abstract}
If a curve in ${\Bbb R^3}$ is closed, then the curvature and the
torsion are periodic functions satisfying some additional constraints.
We show that these constraints can be naturally formulated in terms of
the spectral problem for a $2\times2$ matrix differential operator. This
operator arose in the theory of the self-focusing Nonlinear
Schr\"odinger Equation. 

A simple spectral characterization of Bloch varieties generating 
periodic solutions of the Filament Equation is obtained. We show that
the method of isoperiodic deformations suggested earlier by the
authors for constructing periodic solutions of soliton equations can
be naturally applied to the Filament Equation.
\end{abstract}

\section{Introduction}
\label{sec:introduction}

In our article we shall study the periodic problem for the Filament Equation
\beq
\frac{\partial \vec\gamma(s,t)}{\partial t}=k(s,t) \vec{b}(s,t)
\label{filament_equation}
\eeq
where $\gamma(s,t)$ is $t$-dependent family of smooth curves in 
${\Bbb{R}}^3$, $s$ is a natural parameter on these curves, $k(s,t)$
is the curvature, $\vec{b}(s,t)$ is the binormal vector (the necessary
definitions from differential geometry are collected in 
Section~\ref{sec:curves}).

The Filament Equation describes the motion of a very thin isolated
vortex filament in an incompressible unbounded fluid. It was derived
by Da Rios \cite{DaRios} in the year 1906 and rediscovered in the 60's
by R.~J.~Arms, F.~R.~Hama, R.~J.~Betchov (see the historical article 
\cite{Ric91} for more details). 

By a periodic problem we mean constructing solutions of
(\ref{filament_equation}) such that for any $t=t_0$ the curve 
$\gamma(s,t_0)$ is closed:
\beq
\vec\gamma(s+l,t_0)=\vec\gamma(s,t_0)
\label{periodicity}
\eeq

Without loss of generality we shall assume the length $l$ of the curve
to be equal $2\pi$. A closed curve can be naturally interpreted as a 
smooth isometric map
\beq
\gamma: S^1\rightarrow{\Bbb R}^3
\eeq

In 1972 H.~Hasimoto \cite{Has72} found a change of variables connecting 
(\ref{filament_equation}) with the the self-focusing Nonlinear 
Schr\"odinger Equation (NLS)
\beq
i\frac{\partial q(s,t)}{\partial{t}}+
\frac{\partial^2 q(s,t)}{\partial{s^2}}+ 
\frac12\left|q(s,t)\right|^2q(s,t)=0.
\label{nls}
\eeq

This change of variables associates with a curve $\gamma(s)$ a complex
function $q(s)$:
\beq
{\cal H}: \gamma(s)\rightarrow q(s)=
k(s)e^{i\int^s\kappa(\tilde{s})d\tilde{s}}
\label{hasimoto}
\eeq

(In \cite{DS94} A.~Doliwa and P.~M.~Santini have shown that under some
natural assumptions any integrable motion of a curve in $S^3$ results
in the NLS hierarchy.)

The periodic problem for NLS is well-studied (see Section~\ref{sec:nls}). 
Unfortunately results from the periodic NLS theory can not be applied 
directly to the Filament Equation because the Hasimoto map 
(\ref{hasimoto}) does not map periodic functions to the periodic
ones. It is easy to check, that

\begin{enumerate}
\item For a generic closed curve $\vec\gamma(s+2\pi)=\vec\gamma(s)$
  the corresponding potential $q(s)$ is quasi-periodic:
\beq
q(s+2\pi)=e^{i\phi}q(s).
\eeq
\item For a generic periodic potential $q(s+2\pi)=q(s)$ neither
  $\vec\gamma(s)$ nor the velocity vector $\vec{v}(s)$ are periodic.
\end{enumerate}

The problem of constructing periodic algebro-geometric solutions of
the Filament Equation (\ref{filament_equation}) was studied by 
A.~M.~Calini in her Ph.D. dissertation \cite{Cal94}. In \cite{Cal94} 
a number of interesting results were obtained. In particular explicit exact
solutions were constructed. Unfortunately, no characterization of
periodic solutions of the Filament equation was given in \cite{Cal94}.
As  pointed out by S.~P.~Novikov, without such a characterization the
periodic problem for (\ref{filament_equation}) can not be considered as
completely solved.

Thus we are faced with the following differential-geometrical problem:
If a curve in ${\Bbb R^3}$ is closed, then the curvature and the torsion 
functions are also periodic. But the periodicity of the curvature and the
torsion functions does not imply automatically that the corresponding
curve is closed. To obtain a closed curve we have to add 
additional constraints on the curvature and the torsion. We show that
these constraints can be naturally written in terms of the spectral
problem for a $2\times2$ matrix differential operator associated with
NLS equation (see Sections~\ref{sec:nls},~\ref{sec:periodic}). 

Another interesting question associated with the Hasimoto map arose in
the Hamiltonian theory of the NLS equation. J.~J.~Millson and B.~Zombro
have shown in \cite{MZ96} that the space of smooth isometric maps of
the unit circle into ${\Bbb R}^3$ modulo proper Euclidean motions
has a natural K\"ahler structure. The imaginary part of this structure
coincides with one of the higher NLS symplectic structures. In order to
study this relation it is important to have an explicit description of
the spaces connected by the Hasimoto map. We belive that the characterization
obtained in Section~\ref{sec:periodic} may be applied to this problem.

The problem of characterizing the spectral data corresponding to the
periodic NLS solutions it rather non-trivial. A convenient approach to
this problem based on the so-called isoperiodic deformations was
suggested by the authors in \cite{GS95}. In Section~\ref{sec:periods}
we show that this approach can be naturally applied to the Filament
Equation. 

One of the authors (P.G.) would like to express his gratitude to
Prof. S.~P.~Novikov for the invitation to several visits of the Maryland
University (the last one was in the fall 1996) and for the interest to
this work. He is also grateful to Prof. J.~J.~Millson and B.~Zombro 
for numerous discussions about this problem. 

\section{Curves in 3-dimensional Euclidean space}
\label{sec:curves}
Let $\vec\gamma(s)$ be a smooth, parameterized curve in Euclidean 3-space:
\beq
\vec\gamma(s)=(x^1(s),x^2(s),x^3(s))\in{\Bbb R}^3.
\eeq

Denote by $\vec{v}(s)$ and $\vec{w}(s)$ the velocity and the
acceleration respectively:
\beq
\vec{v}(s)=\frac{d\vec\gamma(s)}{ds}=\left(\frac{dx^1(s)}{ds},\frac{dx^2(s)}
{ds},\frac{dx^3(s)}{ds}\right)
\eeq
\beq
\vec{w}(s)=\frac{d^2\vec\gamma(s)}{ds^2}=
\left(\frac{d^2x^1(s)}{ds^2},\frac{d^2x^2(s)}{ds^2},
\frac{d^2x^3(s)}{ds^2}\right)
\eeq
We shall assume that s is the natural parameter, i.e. the length along
the curve. In other words, the velocity has unit length:
\beq
\left|\vec{v}(s)\right|^2=
\left(\frac{dx(s)}{ds}\right)^2+\left(\frac{dy(s)}{ds}\right)^2+
\left(\frac{dz(s)}{ds}\right)^2=1
\eeq
Then the acceleration vector is orthogonal to the velocity
\beq
<\vec{v}(s),\vec{w}(s)>=0,
\eeq
where
$<\ ,\ >$ denotes the standard scalar product in ${\Bbb R}^3$
\beq
<\vec{v},\vec{w}>=v^1w^1+v^2w^2+v^3w^3.
\eeq
The magnitude of the acceleration vector is called the {\bf curvature}
of the curve
\beq
k(s)=\left|\vec{w}(s)\right|=\sqrt{<\vec{w}(s),\vec{w}(s)>}.
\eeq
For each value of $s$ such that $k(s)\ne0$ we have a natural 
orthogonal reference frame $\left(\vec{v}(s),\vec{n}(s),\vec{b}(s)\right)$, 
where 
\beq
\vec{n}(s)=\frac{\vec{w}(s)}{\left|\vec{w}(s)\right|},\ 
\vec{b}(s)=\vec{v}(s)\times\vec{n}(s).
\eeq
Here $\times$ denotes the vector product in ${\Bbb R}^3$
\beq
\vec{v}\times\vec{n}=\left(v^2n^3-v^3n^2,v^3n^1-v^1n^3,v^1n^2-v^2n^1\right)
\eeq
The vector $\vec{n}(s)$ is called the principal normal to the curve,
the vector $\vec{b}(s)$ is called the binormal.

The natural reference frame $\left(\vec{v}(s),\vec{n}(s),\vec{b}(s)\right)$
satisfies the Serret-Frenet equations (see for example \cite{DFN84})
$$
\frac{d\vec{v}(s)}{ds}=k(s)\vec{n}(s)
$$
\beq
\frac{d\vec{n}(s)}{ds}=-k(s)\vec{v}(s)+\kappa(s)\vec{b}(s)
\label{frenet}
\eeq
$$
\frac{d\vec{b}(s)}{ds}=-\kappa(s)\vec{n}(s)
$$

The function $\kappa(s)$ is called the {\bf torsion}.

Further we shall use another reference frame $\left(\vec{e}_1(s), 
\vec{e}_2(s),\vec{e}_3(s)\right)$ associated with a smooth curve
$\vec\gamma(s)$ parameterized by the natural parameter. Let
\beq
\barr{l}
\vec{e}_1(s)=\vec{v}(s)\\ \\
\vec{e}_2(s)=\cos(\theta(s))\ \vec{n}(s)-\sin(\theta(s))\ \vec{b}(s)\\ \\
\vec{e}_3(s)=\sin(\theta(s))\ \vec{n}(s)+\cos(\theta(s))\ \vec{b}(s)
\earr
\eeq
where 
\beq
\theta(s)=\int^s\kappa(\tilde{s})d\tilde{s}.
\label{theta}
\eeq

This reference frame satisfies

$$
\frac{d\vec{e}_1(s)}{ds}=
k(s)\cos(\theta(s))\ \vec{e_2}+k(s)\sin(\theta(s))\ \vec{e_3}(s) \\ \\
$$
\beq
\frac{d\vec{e_2}(s)}{ds}=-k(s)\cos(\theta(s))\ \vec{e_1}(s) 
\label{frenet_2}
\eeq
$$
\frac{d\vec{e_3}(s)}{ds}=-k(s)\sin(\theta(s))\ \vec{e_1}(s)
$$

Equations (\ref{frenet_2}) are written in terms of the Lie algebra
$so(3)$. Let us rewrite them in terms of $su(2)$.

Let $I$, $J$, $K$ be a basis in the space of skew-hermitian matrices
\beq
I=i\sigma_z=\left[ 
{\begin{array}{cc}
{i} & 0 \\
0 &  - {i}
\end{array}}
 \right],\ \ 
J=-i\sigma_y=\left[ 
{\begin{array}{rr}
0 & -1 \\
1 & 0
\end{array}}
 \right],\ \ 
K=-i\sigma_x=\left[ 
{\begin{array}{cc}
0 &  - {i} \\
 - {i} & 0
\end{array}}
 \right], 
\label{basic_quaternions}
\eeq
where $\vec\sigma=(\sigma_x,\sigma_y,\sigma_z)$ are the Pauli matrices.

These matrices satisfy the multiplication rules for the basic quaternions
\beq
IJ=-JI=K, \ JK=-KJ=I,\ KI=-IK=J,\ I^2=J^2=K^2=-1,
\eeq
and they form an orthonormal reference frame in the space of 
$2\times2$ skew-hermitian matrices with the following scalar product
\beq
<A,B>=-\frac12 \hbox{trace}AB.
\label{scalar_2}
\eeq

It is easy to check that the following map:
\beq
\vec{w}=(w^1,w^2,w^3)\rightarrow\hat{W}=w^1\cdot I+w^2\cdot J+w^3\cdot K.
\label{isomorphism_1}
\eeq
is an isometry between the Euclidean space ${\Bbb R}^3$ and the space
of $2\times2$ skew-hermitian matrices with inverse given by
$$
w^1=-\frac12\hbox{trace}\,\hat{W}I, \ w^2=-\frac12\hbox{trace}\,\hat{W}J,
\ w^3=-\frac12\hbox{trace}\,\hat{W}K.
$$

It generates the famous map 
$SU(2) \rightarrow SO(3)$: any unitary matrix $g\in SU(2)$ generates
an isometric map of the space of $2\times2$ hermitian matrices
\beq
x\rightarrow gxg^{-1}.
\label{isomorphism_2}
\eeq 
It is well-known (see for example \cite{DFN84} p. 431) that any isometric 
orientation preserving map of the the space of $2\times2$ hermitian matrices
is generated by (\ref{isomorphism_2}) and the matrix $g$ is defined
uniquely up to multiplication by $-1$. 

Denote the matrices corresponding to the vectors $\vec e_1(s)$, $\vec e_2(s)$
$\vec{e}_3(s)$ by $\hat{E}_1(s)$, $\hat{E}_2(s)$, $\hat{E}_3(s)$
respectively. They form an orthonormal reference frame. Thus there
exists a $2\times2$ unitary matrix $\Omega(s)$ such that
\beq
\hat{E}_1(s)=\Omega^{-1}(s) I \Omega(s), \ 
\label{E_1}
\eeq
$$
\hat{E}_2(s)=\Omega^{-1}(s) J \Omega(s), \
$$
$$
\hat{E}_3(s)=\Omega^{-1}(s) K \Omega(s).
$$

In terms of these $2\times2$ matrices the system of equations 
(\ref{frenet_2}) takes the form 
\beq
\barr{l}
\left[I,\omega(s)\right]=k(s)\cos(\theta(s))\ J+k(s)\sin(\theta(s))\ K \\ \\
\left[J,\omega(s)\right]=-k(s)\cos(\theta(s))\ I \\ \\
\left[K,\omega(s)\right]=-k(s)\sin(\theta(s))\ I
\earr
\label{frenet_3}
\eeq
where
\beq
\omega(s)=\frac{d\Omega(s)}{ds}\Omega^{-1}(s).
\eeq
Taking into account that $\omega(s)$ is skew-hermitian we obtain
\beq
\omega(s)=-k(s)\cos(\theta(s))\ K+k(s)\sin(\theta(s))\ J=
\left[\barr{cc} 0 & \frac{iq(s)}2 \\ \frac{i\bar{q}(s)}2 & 0 \earr \right]
\label{omega}
\eeq
where 
\beq
q(s)=k(s)e^{i\theta(s)}.
\label{hasimoto_2}
\eeq

Starting from a curve $\vec\gamma(s)\in{\Bbb R}^3$ we have thus constructed
a potential $q(s)$. Formula (\ref{hasimoto_2}) coincides with the 
Hasimoto transformation (\ref{hasimoto}).

Let us discuss the inverse map (the map from the space of complex-valued
functions of one real variable to the space of curves in ${\Bbb R}^3$). 

\begin{lemma}
\label{lemma:reconstruction_1}
Let $q(s)$ be a complex-valued smooth function of one real variable
$s$ such that $q(s)\ne0$ for all $s$. Then there exists an unique (up
to a proper isometry of ${\Bbb R}^3$) curve $\gamma(s)$ such that 
${\cal H}\gamma(s)=e^{i\phi}q(s)$ where ${\cal H}$ is the Hasimoto
map. (Recall that the image of the Hasimoto map is defined up
to an arbitrary constant phase $\phi$). The curve $\gamma(s)$ can be
constructed by using the following procedure:

\begin{enumerate}
\item Define a matrix $\omega(s)$ by (\ref{omega}).
\item \label{step_2}
  Define a $2\times2$ matrix function $\Omega(s)$ as a
  solution of the following linear ordinary differential equation:
\beq
\frac{d}{ds}\Omega(s)=\omega(s)\Omega(s)
\label{lin_problem_1}
\eeq
such that $\Omega(0)$ is an unitary matrix. (From
(\ref{lin_problem_1}) it follows that the matrix $\Omega(s)$ is
unitary for all $s$).
\item Let $\hat{E}_1(s)$ be a skew-hermitian matrix defined by the
  formula (\ref{E_1}), and $\vec{e}_1(s)$ be the corresponding vector 
  in ${\Bbb R}^3$.
\item The curve $\vec\gamma(s)$ and the corresponding matrix-valued
  function $\hat\Gamma(s)$ are defined by:
\beq
\vec\gamma(s)=\vec\gamma(0)+\int_0^s \vec{e}_1(\tilde s)d\tilde s, \
\hat\Gamma(s)=\hat\Gamma(0)+\int_0^s \hat{E}_1(\tilde s)d\tilde s
\label{integral}
\eeq
\end{enumerate}
\end{lemma}

The proof of Lemma~\ref{lemma:reconstruction_1} is standard, 
so we will not present it here.

\begin{remark} 
If the potential $q(s)$is known, so is the curvature and 
the torsion. Indeed one has
\beq
k(s)=|q(s)|, \ \kappa(s)=\frac{d}{ds} \hbox{arg}\,q(s), 
\eeq
and the curve $\gamma(s)$ is defined uniquely up to an
isometry (see for example \cite{DFN84}). However the reconstruction
procedure described above is essential for the approach used in this 
article.
\end{remark}

From Lemma~\ref{lemma:reconstruction_1} it follows directly:
\begin{lemma} 
\label{lemma:periodicity}
The curve $\vec\gamma(s)$ constructed
in Lemma~\ref{lemma:reconstruction_1} is periodic with period
  $l=2\pi$ if and only if the following two conditions are fulfilled:
\begin{enumerate}
\item The matrix $\hat{E}_1(s)$ defined by the formula (\ref{E_1})
  is periodic with period $2\pi$
\beq
\hat{E}_1(s+2\pi)\equiv\hat{E}_1(s).
\label{cond_1}
\eeq
\item The integral of $\hat{E}_1(s)$ over one period vanishes, i.e.
\beq
\int_0^{2\pi} \hat{E}_1(s)ds = \left[\barr{cc}0 & 0 \\ 0 & 0\earr \right].
\label{cond_2}
\eeq
\end{enumerate}
\end{lemma}

\section{Periodic theory of Nonlinear Schr\"odinger Equation}.
\label{sec:nls}

\subsection{The zero-curvature representation}
\label{sec:zero-curvature}

The self-focusing Nonlinear Schr\"odinger Equation (NLS)
\beq
i\frac{\partial q(s,t)}{\partial{t}}+
\frac{\partial^2 q(s,t)}{\partial{s^2}}+ 
\frac12\left|q(s,t)\right|^2q(s,t)=0.
\label{nls_2}
\eeq
is one of the most important soliton equations. In 1971 V.~E.~Zakharov
and A.~B.~Shabat \cite{ZS71} proved that NLS can be integrated by the
inverse scattering method. This method is based on the so-called 
zero-curvature representation for NLS (a good introduction to the NLS
theory can be found in the book \cite{FT80} by L.~D.~Faddeev and
L.~A.~Takhtajian): 

Consider the following pair of linear problems:
\beq
\frac{\partial F(x,t,\lambda)}{\partial{x}}=U(x,t,\lambda)F(x,t,\lambda),
\label{u_problem}
\eeq
\beq
\frac{\partial F(x,t,\lambda)}{\partial{t}}=V(x,t,\lambda)F(x,t,\lambda)
\label{v_problem}
\eeq
where $F(x,t,\lambda)$ is a vector-valued function
\beq
F(x,t,\lambda)=\left(\barr{c} f_1(x,t,\lambda) \\ f_2(x,t,\lambda)
\earr\right),
\eeq
$U(x,t,\lambda)$, $V(x,t,\lambda)$ are the following $2\times2$
matrices, which depend  polynomially on spectral parameter $\lambda$:
\beq
U(x,t,\lambda)=\left[ 
{\begin{array}{cc}
 - \,{\displaystyle \frac {1}{2}}\,{i}\,{ \lambda} & 
{\displaystyle \frac {1}{2}}\,{i}\,q(\,{x}, {t}\,) \\ [2ex]
{\displaystyle \frac {1}{2}}\,{i}\,\bar{q}(\,{x}, {t}\,) & 
{\displaystyle \frac {1}{2}}\,{i}\,{ \lambda}
\end{array}}
 \right] 
\eeq
\beq
V(x,t,\lambda)=\left[ 
{\begin{array}{c}
{\displaystyle \frac {1}{4}}\,{i}\,{q}(\,{x}, {t}\,)\,\bar{q}
(\,{x}, {t}\,) - {\displaystyle \frac {1}{2}}\,{i}\,{ \lambda}^{2
}\,,  - \,{\displaystyle \frac {1}{2}}\, \left( \! \,{\frac {{ 
\partial}}{{ \partial}{x}}}\, q(\,{x}, {t}\,)\, \!  \right) 
 + {\displaystyle \frac {1}{2}}\,{i}\,q(\,{x}, {t}\,)\,{ 
\lambda} \\ [2ex]
{\displaystyle \frac {1}{2}}\, \left( \! \,{\frac {{ \partial}}{{
 \partial}{x}}}\,\bar{q}(\,{x}, {t}\,)\, \!  \right)  + 
{\displaystyle \frac {1}{2}}\,{i}\,\bar{q}(\,{x}, {t}\,)\,{ 
\lambda}\,, \, - \,{\displaystyle \frac {1}{4}}\,{i}\,q(\,{
x}, {t}\,)\,\bar{q}(\,{x}, {t}\,) + {\displaystyle \frac {1}{2}}
\,{i}\,{ \lambda}^{2}
\end{array}}
 \right] 
\eeq
The system (\ref{u_problem}), (\ref{v_problem}) is compatible if and only if
\beq
\frac{\partial U(x,t,\lambda)}{\partial t}- 
\frac{\partial V(x,t,\lambda)}{\partial x}+
\left[U(x,t,\lambda),V(x,t,\lambda)\right]=0.
\label{comatibility}
\eeq
(Here $[\ ,\ ]$ is the standard matrix commutator $[A,B]=AB-BA$.)
A simple direct calculation shows, that (\ref{comatibility}) is
equivalent to (\ref{nls}).

The representation (\ref{comatibility}) is called the {\bf
  zero-curvature} representation (see e.g. \cite{FT80} 
  for additional information).

\begin{remark}
The word ``zero-curvature'' means the following. The matrices 
$U(x,t,\lambda)$, $V(x,t,\lambda)$ can be interpreted as local
connection coefficients in the trivial bundle with base ${\Bbb R}^2$
and fiber ${\Bbb C}^2$. Here $(x,t)$ is a point of the base,
$F(x,t,\lambda)$ takes values in the fiber and $\lambda$ is
a parameter. Then (\ref{comatibility}) means exactly that the
curvature of this connection vanishes for all $\lambda$.
\end{remark}

\subsection{The auxiliary linear problem and gauge transformations}
\label{sec:auxiliary}

The linear problem (\ref{u_problem}) is called the {\bf auxiliary linear
problem} for the NLS equation. It plays a crucial role in the inverse 
scattering method. We have a spectral problem for a first-order ordinary
differential operator in the variable $x$, with a spectral parameter 
$\lambda$. This operator also depend on an additional
parameter $t$. We shall study this spectral problem for a fixed
$t=t_0$. Henceforth we will drop the $t$-dependence from the
notations and write $U(x,\lambda)$, $F(x,\lambda)$,
$q(x)$ instead of $U(x,t_0,\lambda)$, $F(x,t_0,\lambda)$,
$q(x,t_0)$ respectively. Also we will  rewrite (\ref{u_problem}) in the
following form
\beq
L(\lambda)F(x,\lambda)=0,\ \  L(\lambda)=\frac{d}{dx}-\left[ 
{\begin{array}{cc}
 - \,{\displaystyle \frac {1}{2}}\,{i}\,{ \lambda} & 
{\displaystyle \frac {1}{2}}\,{i}\,q(x) \\ [2ex]
{\displaystyle \frac {1}{2}}\,{i}\,\overline{q(x)} & 
{\displaystyle \frac {1}{2}}\,{i}\,{ \lambda}
\end{array}}
 \right] 
\label{L_operator}
\eeq

We are looking for a characterization of potentials corresponding to
periodic curves of length $l=2\pi$. Thus we shall study 
the direct spectral transform for the problem (\ref{L_operator}) 
in the class of smooth complex-valued potentials $q(x)$ such that
\beq
q(s+2\pi)=e^{i\phi}q(s), \ \ \phi\in{\Bbb R}.
\label{per_cond_1}
\eeq
(It is well-known that $\phi$ is an integral of motion for the
filament equations). 

Usually the quasi-periodic spectral theory is much more complicated
than the periodic one. Fortunately in the present case we can handle
the situation (this fact was pointed out in \cite{FT80}). The reason
is that the linear problem (\ref{L_operator}) is invariant under the 
following gauge transformations:
\beq
F(x,\lambda)=\left(\barr{c}f_1(x,\lambda)\\f_2(x,\lambda)\earr\right)
\rightarrow\tilde{F}(x,\tilde\lambda)=
\left(\barr{c}e^{\frac{i}2\alpha{x}}f_1(x,\lambda)\\
e^{-\frac{i}2\alpha{x}}f_2(x,\lambda)\earr\right),
\label{gauge}
\eeq
\beq
\lambda\rightarrow\tilde\lambda=\lambda-\alpha,\
q(x)\rightarrow\tilde{q}(x)=e^{i\alpha{x}}q(x)
\label{gauge_2}
\eeq
\beq
\tilde L(\tilde\lambda)\tilde{F}(x,\tilde\lambda)=0,\ \  
\tilde L(\tilde\lambda)=\frac{d}{dx}-\left[ 
{\begin{array}{cc}
 - \,{\displaystyle \frac {1}{2}}\,{i}\,{ \tilde\lambda} & 
{\displaystyle \frac {1}{2}}\,{i}\,\tilde{q}(x) \\ [2ex]
{\displaystyle \frac {1}{2}}\,{i}\,\overline{\tilde{q}(x)} & 
{\displaystyle \frac {1}{2}}\,{i}\,{ \tilde\lambda}
\end{array}}
 \right] 
\label{L_tilde_operator}
\eeq
parameterized by $\alpha\in{\Bbb R}$.

With the choice
\beq
\alpha=-\frac{\phi}{2\pi}
\eeq
we have transformed (\ref{L_operator}) to a spectral problem with a purely
periodic potential $\tilde q(x+2\pi)=\tilde q(x)$.

The gauge transformation (\ref{gauge}), (\ref{gauge_2}) respects the
formula (\ref{E_1}). More precisely let $\Omega(x)$ be a $2\times2$
matrix solution of the equation
\beq
L(0)\Omega(x)=0.
\eeq
Then 
\beq
\tilde\Omega(x)=e^{-\frac{i\phi}{4\pi}x\,\sigma_z}\Omega(x), \ \ 
\sigma_z=\left[\barr{cc}1&0\\ 0&-1 \earr\right]
\eeq
satisfy
\beq
\tilde L\left(\frac{\phi}{2\pi}\right)\tilde\Omega(x)=0.
\eeq
and
\beq
\hat{E}_1(x)=\Omega^{-1}(x) I \Omega(x)=
\tilde\Omega^{-1}(x) I \tilde\Omega(x).
\label{E_1_tilde}
\eeq

The gauge transformation (\ref{gauge}), (\ref{gauge_2}) shifts the
spectral parameter $\lambda$. Thus the $\lambda=0$ eigenfunctions of
the problem (\ref{L_operator}) with an arbitrary $\phi$ are equivalent
to eigenfunctions of the purely periodic problem with an arbitrary real
$\lambda$. 

Taking into account the gauge transformation properties of 
(\ref{L_operator}) we obtain
the following modification of the reconstruction procedure described
in Lemma~\ref{lemma:reconstruction_1}.

\begin{lemma}
\label{lemma:reconstruction_2}
Let $q(x)$ be a complex-valued smooth function of one real variable
$x$ such that $q(x)\ne0$ for all $x$. Let $\Lambda_0\in{\Bbb R}$ be a real
point in the spectral plane. Then there exists a  curve $\gamma(x)$ 
(unique up to an isometry of ${\Bbb R}^3$) such that 
${\cal H}\gamma(x)=e^{i\phi-i\Lambda_0x}q(x)$ where ${\cal H}$ is the Hasimoto
map. The curve $\gamma(x)$ can be obtained by the following procedure:

\begin{enumerate}
\item Let $\Omega(x)$ be an arbitrary $2\times2$ matrix
  solution of the following equation:
\beq
L(\Lambda_0)\Omega(x)=0
\label{lin_problem_2}
\eeq
such that $\Omega(0)$ is a unitary matrix.
\item Let $\hat{E}_1(x)$ be the skew-hermitian matrix defined by
 \beq
 \hat{E}_1(x)=\Omega^{-1}(x) I \Omega(x)
 \label{E_1_tilde_2}.
 \eeq
\item The function $\hat\Gamma(x)$ is defined by
\beq
\hat\Gamma(x)=\hat\Gamma(0)+\int_0^x \hat{E}_1(\tilde x)d\tilde x.
\label{integral_2}
\eeq
\end{enumerate}
\end{lemma}

\subsection{Periodic spectral problem and Bloch variety}
\label{sec:bloch}

For the remainder of this section we shall assume that $q(x)$ is a
smooth complex-valued periodic function with period $2\pi$.

If we impose no boundary conditions, then the equation (\ref{L_operator}) has
a two-dimensional space of solutions for any complex $\lambda$. 
A point $\lambda\in{\Bbb C}$ belongs to the spectrum of
(\ref{L_operator}) if and only if this space contains at least one
function bounded on the whole $x$-line.

The structure of the spectrum of the problem (\ref{L_operator}) may be
rather complicated (this structure was studied by Y.~Li and
D.~W.McLaughlin in \cite{LM94}). Fortunately we do not have to know this
structure in detail to construct explicit solutions.

Let us fix a basis of solutions $\vec\varphi^{(1)}(x,\lambda)$, 
$\vec\varphi^{(2)}(x,\lambda)$:
\beq
L(\lambda)\vec\varphi^{(1)}(x,\lambda)=
L(\lambda)\vec\varphi^{(2)}(x,\lambda)=0, \ \ 
\vec\varphi^{(1)}(0,\lambda)=\left(\barr{c}1 \\ 0  \earr\right), \ 
\vec\varphi^{(2)}(0,\lambda)=\left(\barr{c}0 \\ 1  \earr\right).
\eeq
Denote by $\Phi(x,\lambda)$ the $2\times2$ fundamental solution of
(\ref{L_operator}):
\beq
\Phi(x,\lambda)=\left[\vec\varphi^{(1)}(0,\lambda)
\vec\varphi^{(2)}(0,\lambda)\right].
\eeq
The operator $L(\lambda)$ commutes with the shift operator
\beq
f(x)\rightarrow f(x+2\pi).
\eeq
Thus we have
\beq
\Phi(x+2\pi,\lambda)=\Phi(x,\lambda)T(\lambda), \ \hbox{where} \  
T(\lambda)=\Phi(2\pi,\lambda)
\eeq
The matrix $T(\lambda)$ is called the monodromy matrix. For generic
$\lambda$ it can be diagonalized. The common eigenfunctions of
$L(\lambda)$ and of the shift operator 
$$
\vec\psi^{(1)}(x,\lambda)=
\left(\barr{c}\psi^{(1)}_1(x,\lambda)\\ \psi^{(1)}_2(x,\lambda)\earr\right), 
\ \ \ \
\vec\psi^{(2)}(x,\lambda)=
\left(\barr{c}\psi^{(2)}_1(x,\lambda)\\ \psi^{(2)}_2(x,\lambda)\earr\right)
$$
\beq
L(\lambda)\vec\psi^{(1)}(x,\lambda)=0, \ \ 
\vec\psi^{(1)}(x+2\pi,\lambda)=w^{(1)}(\lambda)
\vec\psi^{(1)}(x,\lambda), 
\eeq
$$
L(\lambda)\vec\psi^{(2)}(x,\lambda)=0, \ \ 
\vec\psi^{(2)}(x+2\pi,\lambda)=w^{(2)}(\lambda)
\vec\psi^{(2)}(x,\lambda)
$$
are called the {\bf Bloch functions}. The functions $p^{(1)}(\lambda)$,
$p^{(2)}(\lambda)$ are called quasimomenta. The matrix
$U(x,\lambda)$ is traceless, thus
$$
w^{(1)}(\lambda)w^{(2)}(\lambda)=1.
$$

A point $\lambda\in{\Bbb C}$ belongs to the spectrum of
(\ref{L_operator}) if and only if 
$$
\left|w^{(1)}(\lambda)\right|=1.
$$

It is convenient to fix a normalization of the Bloch functions by
the condition
\beq
\psi^{(1)}_1(0,\lambda)+\psi^{(1)}_2(0,\lambda)=
\psi^{(2)}_1(0,\lambda)+\psi^{(2)}_2(0,\lambda)=1.
\label{bloch_normalization}
\eeq

It is easy to check that
\begin{enumerate}
\item The matrix $T(\lambda)$ is holomorphic in $\lambda$ in the whole  
$\lambda$-plane.
\item $\det T(\lambda)\equiv1$.
\item For all $\lambda\in{\Bbb R}$ $T(\lambda)$ is unitary, i.e. 
$T^{-1}(\lambda)=T^*(\lambda)$ where $*$ denotes hermitian
conjugation. 
\end{enumerate}

For a generic $\lambda\in{\Bbb C}$ $T(\lambda)$ has two eigenvalues 
$w^{(1)}(\lambda)$, $w^{(2)}(\lambda)$ and a pair of corresponding
Bloch functions. In fact there is a holomorphic function $w(\mu)$ on a
hyperelliptic Riemann surface $Y$ and $w^{(1)}(\lambda)=w(\mu_1)$, 
$w^{(2)}(\lambda)=w(\mu_2)$ where $\mu_1$ and $\mu_2$ are the
pre-images of the point $\lambda$ under the projection
$Y\rightarrow{\Bbb C}$ (a Riemann surface is 
called hyperelliptic if is a two-sheeted ramified covering of the
Riemann sphere). For generic potentials the surface $Y$ is connected,
but there exist exceptional potentials such that $Y={\Bbb C}\cup{\Bbb C}$. 

Any hyperelliptic Riemann surface has a natural holomorphic involution
given transposing the sheets. Let us denote this involution by $\sigma$:
\beq
\sigma\lambda=\lambda, \ \ \sigma{w}(\mu)=w^{-1}(\mu).
\eeq

Let 
\beq
\vec\psi(x,\mu)=
\left(\barr{c}\psi_1(x,\mu)\\ \psi_2(x,\mu)\earr\right)
\eeq
be a Bloch solution of $L(\lambda(\mu))\vec\psi(x,\mu)=0$, 
$\vec\psi(x+2\pi,\mu)=w(\mu)\vec\psi(x,\mu)$.
Then 
\beq
\vec\psi^+(x,\mu)=
\left(\barr{c}-\bar\psi_2(x,\mu)\\ \bar\psi_1(x,\mu)\earr\right)
\eeq
is a Bloch solution of $L(\bar\lambda(\mu))\vec\psi^+(x,\mu)=0$, 
$\vec\psi^+(x+2\pi,\mu)=\bar{w}(\mu)\vec\psi^+(x,\mu)$.

Thus on $Y$ there is in addition an antiholomorphic involution
$\sigma\tau$ (we use this notation for historical reasons):
\beq
\sigma\tau\lambda=\bar\lambda, \ \ \sigma\tau{w}(\mu)=\bar{w}(\mu).
\eeq

The Bloch solution of (\ref{L_operator}) $\vec\psi(x,\mu)$ with the
normalization (\ref{bloch_normalization}) is meromorphic in $\mu$ on 
$Y$ with poles which do not depend on $x$.

The function 
\beq
p(\mu)=\frac{1}{2\pi i} \ln w(\mu)
\eeq
is called the {\bf quasimomentum function}. It's differential
\beq
dp(\mu)=\left[\frac{d}{d\mu} p(\mu)\right] d\mu
\eeq
is called the {\bf quasimomentum differential}. Of course $p(\mu)$ is
defined up to adding an arbitrary integer. For a generic potential
$q(x)$ the function 
$p(\mu)$ is essentially multivalued on $Y$, i.e. it's increment is 
non-zero along some cycles in $Y$. The quasimomentum differential is a
well-defined holomorphic differential on the finite part of $Y$. 

The Riemann surface $Y$ is called the {\bf Bloch variety}. It plays a
crucial role in the inverse problem for periodic potentials.
The structure of $Y$ was studied in details by one of the authors (M.S.) in
\cite{Sch96}. Let us recall some basic facts.

A point $\lambda\in{\Bbb C}$ is called regular if 
$w^{(1)}(\lambda)\ne w^{(2)}(\lambda)$ and irregular if 
$w^{(1)}(\lambda)= w^{(2)}(\lambda)$. 
We shall distinguish 3 types of irregular points.
\begin{enumerate}
\item Branch points. 
\item Non-removable double points.
\item Removable double points.
\end{enumerate}

An irregular point $\lambda_0$ is called a branch point if going
around this point we come from one sheet of $Y$ to the other one (i.e. the
monodromy around this point is non-trivial). If the monodromy around
an irregular point $\lambda_0$ is trivial, then $\lambda_0$ is called
a double point. In a neighbourhood of a double point we have a pair 
of locally holomorphic functions $w^{(1)}(\lambda)$, $w^{(2)}(\lambda)$,  
and the Bloch functions $\vec\psi^{(1)}(x,\lambda)$,
$\vec\psi^{(2)}(x,\lambda)$ are locally meromorphic. A double point
$\lambda_0$ is called non-removable if
$\vec\psi^{(1)}(x,\lambda_0)=\vec\psi^{(2)}(x,\lambda_0)$, and
removable otherwise.  

\begin{lemma}
\label{lemma:real_lambda}
Let $\lambda$ be a real point of the spectral plane, i.e.  
$\bar\lambda=\lambda$.Then:
\begin{enumerate}
\item $\hbox{tr}\, T(\lambda)$ is a real function of $\lambda$ and 
$\left|\hbox{tr}\, T(\lambda)\right|\le2$.
\item $\lambda$ lies in the spectrum of (\ref{L_operator}).
\item The point $\lambda$ is regular if and only if 
$\left|\hbox{tr}\, T(\lambda)\right|<2$. 
\item If $\hbox{tr}\, T(\lambda)=\pm2$ then the matrix $T(\lambda)$ is diagonal
\beq
T(\lambda)=\pm \left[ {\begin{array}{cc} {1} & 0 \\
0 & {1} \end{array}} \right]
\eeq
and $\lambda$ is a removable double point.
\end{enumerate}
\end{lemma}
The fact that all irregular points on the real line are removable
double points was proved in \cite{Sch96}. All other statements of
Lemma~\ref{lemma:real_lambda} follow directly from the unitarity 
of $T(\lambda)$.

For generic $q(x)$ $Y$ has infinitely many branch points (and the
genus of $Y$ is infinite). But the asymptotic structure of $Y$ near
infinity is rather simple. Let $\varepsilon$ be a real positive
constant. Then there exists a constant $R$ depending on $q(x)$ and
$\varepsilon$ such that:
\begin{enumerate}
\item If $\lambda\in{\Bbb C}$, $|\lambda|>R$ and
  $|\lambda-k|>\varepsilon$ for any integer $k\in{\Bbb Z}$, then the point
  $\lambda$ is regular.
\item Let $k\in{\Bbb Z}$ be an integer such, that $|k|>R$. Then the
  $\varepsilon$-neighbourhood of the point $k$ in the complex plane
  contains either a pair of complex-conjugate  branch points or one
  removable real double point.
\end{enumerate}

A potential $q(x)$ is called {\bf finite-gap} or {\bf
  algebro-geometric} if $Y$ has only a finite number
of branch points (the number of non-removable double point is always
finite). A finite-gap potential has an explicit representation in terms of
the Riemann $\theta$-functions. Using the methods of \cite{Sch96} and 
\cite{GS95} it is rather easy to prove that any smooth potential can
be approximated by the finite-gap potentials. 

Finite-gap potentials in the context of the soliton theory first appeared
in the soliton theory in the article \cite{Nov74} by S.~P.~Novikov (see
the book \cite{ZMNP} for additional information). The finite-gap
theory of first-order matrix differential operators, including the
$\theta$-function formulas was developed by B.~A.~Dubrovin (see
\cite{Dub77}). The algebro-geometrical solutions of the NLS equation
were studied by E.~Previato in \cite{Pre85}.

An analog of the finite-gap theory for generic periodic potentials can
also be developed. For the first-order matrix operators including the
NLS case it was done by one of the authors (M.S.) in \cite{Sch96}.

Let us recall some useful formulas from \cite{Sch96}.

\begin{lemma}
\label{lemma:quasimomentum}
\begin{enumerate}
\item 
  \label{quasi:part_1}
  Let $\lambda\in{\Bbb C}$ be a regular or a removable double
  point and let $\mu_1$ and $\mu_2=\sigma\mu_1$ be the pre-images of
  $\lambda$. Then
\beq
dp(\mu_1)=-\frac
{\int\limits_0^{2\pi}
\left[\psi_1(x,\mu_1)\psi_2(x,\mu_2)+\psi_2(x,\mu_1)\psi_1(x,\mu_2)
\right] dx}{4\pi\left[
\psi_1(0,\mu_1)\psi_2(0,\mu_2)-\psi_2(0,\mu_1)\psi_1(0,\mu_2)\right]}
d\lambda(\mu_1)
\label{dp_1}
\eeq
\item 
  \label{quasi:part_2}
  Let $\lambda\in{\Bbb C}$ be a removable double point, $\mu_1$
  be one of the pre-images of $\lambda$. Then 
\beq
\int_0^{2\pi} \psi_1(x,\mu_1)\psi_2(x,\mu_1) dx =0.
\label{dp_2}
\eeq
\end{enumerate}
\end{lemma}
Proof of Lemma~\ref{lemma:quasimomentum}. 
Let $\mu$, $\nu$ be an arbitrary pair of points in $Y$ and
$\lambda(\mu)$, $\lambda(\nu)$ their projections to the
$\lambda$-plane. A direct calculation shows that
$$
\frac{d}{dx}
\left[\psi_1(x,\mu)\psi_2(x,\nu)-\psi_2(x,\mu)\psi_1(x,\nu)
\right]=
$$
\beq
=-\frac{i}{2}(\lambda(\mu)-\lambda(\nu))
\left[\psi_1(x,\mu)\psi_2(x,\nu)+\psi_2(x,\mu)\psi_1(x,\nu)
\right].
\label{quasi_1}
\eeq
Integrating (\ref{quasi_1}) over a period we get
$$
\left(w(\mu)w(\nu)-1\right)
\left[\psi_1(0,\mu)\psi_2(0,\nu)-\psi_2(0,\mu)\psi_1(0,\nu)
\right]=
$$
\beq
-\frac{i}{2}(\lambda(\mu)-\lambda(\nu))\int_0^{2\pi}
\left[\psi_1(x,\mu)\psi_2(x,\nu)+\psi_2(x,\mu)\psi_1(x,\nu)
\right]
\label{quasi_2}
\eeq
To prove part~\ref{quasi:part_1} let us assume that $\nu=\mu_2$,
$\mu=\mu_1+\delta$, $\delta\rightarrow0$. 
From the fact $\lambda$ is a regular or a removable double point it
follows that the Wronskian in the denominator of (\ref{dp_1}) is non-zero.
Then (\ref{dp_1}) follows directly from (\ref{quasi_2}). 

To prove part~\ref{quasi:part_2} let us assume that
$\mu=\mu_1$, $\nu=\mu_1+\delta$ where $\delta\rightarrow0$. The
left-hand side of (\ref{quasi_2}) has at least a second-order zero as 
$\delta\rightarrow0$, $(\lambda(\mu)-\lambda(\nu))$ has a first-order
zero, thus the integral at the right-hand side of (\ref{quasi_2})
vanishes as $\delta\rightarrow0$. This completes the proof.

\section{Riemann surfaces corresponding to periodic curves}
\label{sec:periodic}

We now state the main result of our article.

\begin{theorem} 
\label{main_theorem}
Let $q(x)$ be a complex-valued smooth periodic function of one
variable $x$, $q(x)\ne0$ for all $x$, $q(x+2\pi)=q(x)$,
$\Lambda_0\in{\Bbb R}$ an arbitrary real number, $\hat\Gamma(x)$
the corresponding curve constructed in
Lemma~\ref{lemma:reconstruction_2}. Then
\begin{enumerate}
\item The matrix $\hat E_1(x)$ is periodic with period $2\pi$, i.e.
  $\hat E_1(x+2\pi)=\hat E_1(x)$, if and only if $\Lambda_0$ is a
  double point of the Bloch variety $Y$ (let us recall that any real
  double point is automatically removable).
\item The function $\hat\Gamma(x)$ is periodic with period $2\pi$
  $\hat\Gamma(x+2\pi)=\hat\Gamma(x)$, if and only if $\Lambda_0$ is a
  double point of $Y$ and $dp(\mu_1)=0$, $dp(\mu_2)=0$  where $\mu_1$,
  $\mu_2$ are the 
  pre-images of $\Lambda_0$ under the projection $Y\rightarrow{\Bbb C}$.
\end{enumerate}
\end{theorem}

Proof of Theorem~\ref{main_theorem}. Let $\mu_0$ be one of the
pre-images of $\Lambda_0$, let
\beq
\vec\psi(x,\mu_0)=\left(\barr{c}\psi_1(x,\mu_0) \\ \psi_2(x,\mu_0)
  \earr \right) 
\eeq
be a Bloch solution of (\ref{L_operator}) with a normalization such
that 
\beq
\psi_1(0,\mu_0)\bar\psi_1(0,\mu_0)+\psi_2(0,\mu_0)\bar\psi_2(0,\mu_0)=1.
\eeq
$\Lambda_0$ is real thus 
\beq
\vec\psi(x,\sigma\mu_0)=\left(\barr{c}-\bar\psi_2(x,\mu_0) \\ 
\bar\psi_1(x,\mu_0) \earr \right) 
\eeq
and the matrix 
\beq
\Omega(x)=\left[\barr{cc}\psi_1(x,\mu_0) & -\bar\psi_2(x,\mu_0) \\
\psi_2(x,\mu_0) & \bar\psi_1(x,\mu_0) \earr \right]
\eeq
is unitary, satisfies (\ref{lin_problem_2}) and one has
\beq
\Omega(x+2\pi)=\Omega(x)\left[\barr{cc} w(\mu_0) & 0 \\
0 & w^{-1}(\mu_0) \earr \right].
\eeq

Thus for the function $\hat E_1(x)$ defined by (\ref{E_1_tilde_2}) we
have 
\beq
\hat E_1(x+2\pi)=\left[\barr{cc} w^{-1}(\mu_0) & 0 \\
0 & w(\mu_0) \earr \right] \hat E_1(x)
\left[\barr{cc} w(\mu_0) & 0 \\
0 & w^{-1}(\mu_0) \earr \right].
\label{E_1_Bloch}
\eeq

From (\ref{E_1_Bloch}) it follows that $\hat E_1(x+2\pi)\equiv\hat E_1(x)$ 
if and only if one of the following two conditions is fulfilled:

\begin{enumerate}
\item $\hat E_1(x)$ is diagonal for all $x$.
\item $w(\mu_0)=\pm1$.
\end{enumerate}

If the matrix $\hat E_1(x)$ is diagonal for all $x$ then
$q(x)\equiv0$, but $q(x)$ is everywhere non-zero by assumption. Thus 
$w(\mu_0)=\pm1$ and $\Lambda_0$ is an irregular point. By
Lemma~\ref{lemma:real_lambda} $\Lambda_0$ is therefore a removable
double point.

A simple calculation shows that
\beq
\int_0^{2\pi}\hat E_1(x) dx =\left[\barr{cc} a & -\bar b \\
b & -a \earr \right] 
\eeq
where
\beq
a=i\int_0^{2\pi}\psi_1(x,\mu_0)\bar\psi_1(x,\mu_0)-
\psi_2(x,\mu_0)\bar\psi_2(x,\mu_0)dx,
\eeq
\beq
b=-2i\int_0^{2\pi}\psi_1(x,\mu_0)\psi_2(x,\mu_0)dx.
\eeq

By Lemma~\ref{lemma:quasimomentum}
\beq
b=0, \ \ 
a=-4\pi\left.\frac{\partial p(\mu)}{\partial\lambda(\mu)}\right|_{\mu=\mu_0}.
\eeq
Thus $a=0$ and $\hat\Gamma(x+2\pi)=\hat\Gamma(x)$ if and only if
$dp(\mu_0)=0$. Theorem~\ref{main_theorem} is proved.

\section{Deformations of Bloch varieties and periodic solutions of
  the Filament Equation}
\label{sec:periods}

It is well-known that if the Bloch variety is algebraic, then the 
solutions of the soliton equation can be written
explicitly in terms of the Riemann $\theta$-functions. Such solutions
are called algebro-geometric or finite-gap. But if we start
from a generic algebraic Riemann surface, then the corresponding
$\theta$-functional solutions are quasi-periodic in space. Riemann
surfaces generating purely periodic solutions form a rather complicated
transcendental subvariety in the moduli space of all Riemann surfaces
(let us call it the ``Periodic subvariety''). 

A characterization of this subvariety for the periodic
Korteveg-de~Vries equation in terms of conformal maps was obtained by
V.~A.~Marchenko and I.~V.~Ostrovski in \cite{MO75}. Another approach
to construct periodic solutions of soliton equations is based on the
so-called period preserving deformations of Riemann surfaces and it was
suggested by the authors in \cite{GS95}. Let us recall in brief the
results of \cite{GS95} concerning the self-focusing NLS equation.
To explain our ideas, for the sake of transparency we shall
consider generic algebro-geometric potentials only. The case of general
periodic potentials can also be studied by the same method, but it 
requires some more details from the soliton theory.

Let $q(x)$ be a generic non-singular algebro-geometric NLS potential with
period $2\pi$: $q(x+2\pi)=q(x)$. Then the associated Bloch variety $Y$
has the following properties:
\begin{enumerate}
\item $Y$ has a finite even number of branch point $\lambda_1$,\ldots
  ,$\lambda_{2g+2}$, 
  $\hbox{Im}\, \lambda_k\ne 0$ for all $k$, the enumeration can be
  chosen such that $\lambda_{2k+2}=\bar\lambda_{2k+1}$. $g$ is called
  the genus of $Y$.
\item $Y$ has two distinct points over $\lambda=\infty$. Denote these
  points by $\infty_+$ and $\infty_-$. A local parameter in
  a neighborhoods of these points can be chosen as $\nu=1/\lambda$.
\item $Y$ has no non-removable double points.
\item \label{dp_definition} The quasimomentum differential has the
  following representation:
\beq
\barr{c}
dp(\mu)=-\frac12\frac{q(\lambda)}
{\sqrt{(\lambda-\lambda_1)\ldots(\lambda-\lambda_{2g+2})}}d\lambda, 
\\ \\ \hbox{where} \ \lambda=\lambda(\mu),
\ q(\lambda)=\lambda^{g+1}+q_g\lambda^g+\ldots+q_1\lambda+q_0,
\earr
\label{dp_representation}
\eeq
where the real constants $q_0$, $q_1$, \ldots, $q_g$ are uniquely defined
by:
\beq
\hbox{Im}\,\oint_cdp(\mu)=0, \ \ \hbox{res}\left.\vphantom{\int}\right|_
{\mu=\pm\infty} dp(\mu)=0,
\label{dp_normalisation}
\eeq
where $c$ being an arbitrary closed cycle in $Y$.
\item Let $c$ be an arbitrary closed cycle in $Y$. Then 
\beq
\oint_cdp(\mu)\in {\Bbb Z}.
\label{dp_cond_1}
\eeq
\item Let $c$ be a path in $Y$ connecting the points $\infty_+$ and
  $\infty_-$, let $p(\mu)$ be a primitive of $dp(\mu)$ defined on $c$ with
  the following normalization:

\beq
p(\mu)=-\frac{\lambda(\mu)}2 +O\left(\frac1\lambda\right)\ \hbox{as} \ 
\mu\rightarrow\infty_+. 
\label{dp_cond_2}
\eeq
Then 
\beq
p(\mu)=\frac{\lambda(\mu)}2 +k+O\left(\frac1\lambda\right)\ \hbox{as}
\ \mu\rightarrow\infty_-, \hbox{where} \ k\in {\Bbb Z}.
\label{dp_cond_3}
\eeq
\end{enumerate}

Let $Y$ be an arbitrary hyperelliptic Riemann surface with the
properties 1)-3) listed above. Then we can construct a family of NLS
solutions corresponding to $Y$ (see for example \cite{Pre85}). The real
dimension of this family is $g+1$. Some of these solutions may have
singularities, but it is not important in the present context. The 
quasiperiods of these solutions depend only on $Y$ and are the same
for the whole family. It is well-know that Riemann surfaces generating
purely periodic in $x$-space solutions can be characterized in terms 
of the quasimomentum differential (see for example \cite{GS95}):
\begin{lemma}
\label{lemma:periodic_solutions}
Let $Y$ be a hyperelliptic Riemann surface with the properties 1)-3)
listed above. Let $q(x)$ be one of the potentials corresponding to
$Y$. Then $q(x)$ is periodic with a period $2\pi$ $q(x+2\pi)=q(x)$ if
and only if a meromorphic differential $dp(\mu)$ uniquely defined 
by the formulas (\ref{dp_representation}), (\ref{dp_normalisation})
has the properties (\ref{dp_cond_1})-(\ref{dp_cond_3}). 
\end{lemma}

In the finite-gap theory the differential $dp(\mu)$ is called the
quasimomentum differential for both the periodic and the quasiperiodic
solutions.

Assume now that we have a family of Riemann surfaces $Y(\xi)$ such
that the following equations, introduced by I.M.~Krichever \cite{Kr92},
N.M.~Ercolani, M.G.~Forest, D.W.~McLaughlin, A.~Sinha \cite{EFMS},
and one of the authors (M.S.) \cite{Sch96} are fulfilled:
\beq
\left.\frac{\partial p(\mu,\xi)}{\partial\xi}\right|_{\mu=const}=
\frac{\omega(\mu,\xi)}{d\lambda(\mu)}
\label{deformation_1}
\eeq
where any fixed $\xi$ $\omega(\mu,\xi)=\tilde\omega(\mu,\xi)d\lambda(\mu)$ 
is a meromorphic differential on $Y(\xi)$ having no poles
outside the points $\infty_+$, $\infty_-$ and at most first
order-poles at the infinite points. Let us recall that any such
meromorphic differential can be written in the following form
\beq
\barr{c}
\omega(\mu,\xi)=-\frac12\frac{o(\lambda,\xi)}
{\sqrt{(\lambda-\lambda_1(\xi))\ldots(\lambda-\lambda_{2g+2}(\xi))}}
d\lambda, 
\\ \\ \hbox{where} \ \lambda=\lambda(\mu),
\ o(\lambda,\xi)=o_g(\xi)\lambda^g+\ldots+o_1(\xi)\lambda+o_0(\xi).
\earr
\label{holomorphic_differentials}
\eeq
(We shall assume that $o_k(\xi)$ are smooth functions of
$\xi$). Assume that $Y(0)$ generates periodic potentials, and for all
$\xi$ the surface $Y(\xi)$ admits the antiholomorphic involution
$\sigma\tau$. Then $Y(\xi)$ generates periodic potentials for all $\xi$. 

The proof of this fact is based on 
Lemma~\ref{lemma:periodic_solutions} and the following
observation: The function on the right-hand side of
(\ref{deformation_1}) is single-valued on $Y$ and decays as
$\mu\rightarrow\infty_\pm$. Thus the deformation (\ref{deformation_1})
affects neither the periods of $dp(\mu)$ nor the asymptotics of
$p(\mu)$ at infinity.

It is well-known that equation (\ref{deformation_1}) is equivalent
to the following system of ordinary differential equations on the
branch points (see \cite{FFM})
\beq
\frac{\partial\lambda_k(\xi)}{\partial\xi}=
-\frac{o(\lambda,\xi)}{q(\lambda,\xi)}.
\label{deformation_2}
\eeq
The right-hand side of (\ref{deformation_2}) involves hyperelliptic
functions of the branch points because we have to calculate 
integrals (\ref{dp_normalisation}) over periods of $Y$ to determine 
the coefficients of $q(\lambda,\xi)$. In \cite{GS95} it was shown
that the system (\ref{deformation_2}) can be naturally embedded into a
slightly bigger system of ordinary differential equations with a
rational right-hand side. 

Denote by $\alpha_k(\xi)$, $k=1,\dots,g+1$ the zeroes of the polynomial
$q(\lambda,\xi)$ (we shall call them the zeroes of the quasimomentum
differential). In fact $dp(\mu,\xi)$ has zeroes at both pre-images of
each $\alpha_k(\xi)$. Since we consider only generic surfaces, we may
assume that all $\alpha_k(\xi)$ are pairwise distinct. 

\begin{lemma}
\label{rational_form}
Let $\omega_k(\mu,\xi)$, $k=1,\ldots,g+1$ be the following basis of
differentials on $Y(\xi)$
\beq
\omega_k(\mu,\xi)=\frac{1}{\lambda-\alpha_k(\xi)}dp(\mu,\xi), \ 
\lambda=\lambda(\mu).
\eeq
Let
\beq
\omega(\mu,\xi)=\sum_{k=1}^{g+1} c_k(\xi)\omega_k(\mu,\xi)
\label{omega_definition}
\eeq
where $c_k(\xi)$ are arbitrary complex constants. Then the flow
(\ref{deformation_1}), generated by the differential $\omega(\mu,\xi)$
has the following representation:
\beq
\barr{c} 
\frac{\partial \lambda_{j}(\xi)}{\partial \xi} 
=-\sum\limits_{k=1}^{g+1}\frac{c_{k}(\xi)}
{\lambda_{j}(\xi)-\alpha_{i}(\xi)}
\\
\frac{\partial \alpha_{j}(\xi)}{\partial \xi} 
=\sum\limits_{{k=1,\ldots,g+1}\atop{k\ne\j}}
\frac{c_k(\xi)+c_j(\xi)}{\alpha_{k}(\xi)-\alpha_{j}(\xi)} 
-\frac{1}{2} \sum\limits_{j=1}^{2g+2} 
\frac{c_{k}(\xi)}{\lambda_{j}(\xi)-\alpha_{k}(\xi)} 
\earr
\label{deformation_3}
\eeq 
\end{lemma}
Lemma~\ref{rational_form} was proved in \cite{GS95}. 

For generic complex $c_k(\xi)$ the flow (\ref{deformation_3}) does not
respect the symmetry $\lambda_{2j+2}(\xi)=\bar\lambda_{2j+1}(\xi)$ . 
To construct NLS solution we have to add additional restrictions on
$c_k(\xi)$.  

If $Y(\xi)$ has the symmetry $\lambda_{2j+2}(\xi)=\bar\lambda_{2j+1}(\xi)$, 
then the polynomial $q(\lambda)$ has
real coefficients, and there are two types of zeroes $\alpha_k(\xi)$:
\begin{enumerate}
\item Real zeroes $\alpha_k(\xi)=\bar\alpha_k(\xi)$.
\item Pairs of complex conjugate zeroes $\alpha_l(\xi)=\bar\alpha_k(\xi)$.
\end{enumerate}

\begin{lemma}
\label{real_flows}
\begin{enumerate}
\item Let the functions $c_k(\xi)$ be chosen such that for all $\xi$
  the following conditions are fulfilled
\begin{enumerate}
\item If $\alpha_k(\xi)$ is a real zero of the quasimomentum differential,
  then $c_k(\xi)=\bar{c}_k(\xi)$.
\item If $\alpha_k(\xi)$, $\alpha_l(\xi)$ is a pair of complex
  conjugate zeroes 
  of the quasimomentum differential, then $c_l(\xi)=\bar{c}_k(\xi)$. 
\end{enumerate}
Then the flow (\ref{deformation_3}) respects the symmetry
  $\lambda_{2j+2}(\xi)=\bar\lambda_{2j+1}(\xi)$.
\item Consider the subvariety of all hyperelliptic Riemann surfaces
  generating periodic solutions of the self-focusing NLS in the moduli
  space. Let $Y$ be a generic point of this subvariety. Then the flows
  described above act locally transitive at this subvariety. 
\end{enumerate}
\end{lemma}
The first statement of Lemma~\ref{real_flows} follows directly
from the formulas (\ref{deformation_3}).
To prove the last statement it is sufficient to calculate the number
of conditions in a generic point (see \cite{GS95} for details).
 
The algorithm for constructing periodic solutions of the self-focusing NLS 
suggested in
\cite{GS95} is the following: Assume that we know at least one Riemann
surface with $2g+2$ branch points generating such solutions. Such a
surface can be constructed in a neighbourhood of a constant
solution by methods of perturbation theory. Then integrating the
system of ordinary differential equations (\ref{deformation_3}) we can
construct all NLS solutions in a neighbourhood of the original one. 
(In fact, any NLS solution can be constructed by this method, but to
reach some of them we have to pass through singular points of the
system (\ref{deformation_3}).)

Let us show that the method of \cite{GS95} can be naturally applied to
construct periodic solutions in the $x$-space of the Filament equation. (Let us
recall that if the curve $\vec\gamma(x,t)$ is periodic in $x$ for a
fixed $t=t_0$, then it is periodic in $x$ for all $t$ and the
$x$-period does not depend on $t$.) 

\begin{theorem}
\label{filaments_deformations}
Let $Y$ be a hyperelliptic Riemann surface such that:
\begin{enumerate}
\item $Y$ is a Bloch variety corresponding to a smooth periodic
  potential $q(x)$, $q(x+2\pi)=q(x)$.
\item $Y$ has $2g+2$ branch points and no non-removable double points.
\item There exists a point $\Lambda_0\in{\Bbb R}$ such that
\begin{enumerate}
\item $\Lambda_0$ is a removable double point, i.e. the values of the
  quasimomentum function $p(\mu)$ at the pre-images of $\Lambda_0$
  are integer or half-integer.
\item $\Lambda_0$ coincides with one of the zeroes of the quasimomentum
  differential. Since the zeroes of the quasimomentum differential
  $\alpha_1$, \ldots, $\alpha_{g+1}$ have no natural order,
  without loss of generality one may assume therefore 
  $\Lambda_0=\alpha_{g+1}$.
\end{enumerate}
\end{enumerate}
(Let us recall that by Theorem~\ref{main_theorem} the pair
$(Y,\Lambda_0)$ generates periodic solutions in $x$-space of 
the Filament equation). 

Let $c_k(\xi)$, $k=1,\ldots,g+1$ be arbitrary smooth
complex-valued functions of a real variable $\xi$ defined in a
neighbourhood of the point $\xi=0$ such, that: 
\begin{enumerate}
\item For all $\xi$ the $c_k(\xi)$ satisfy the reality conditions of
  Lemma~\ref{real_flows}.
\item $c_{g+1}(\xi)\equiv0$. 
\end{enumerate}

Let $\lambda_k(\xi)$, $k=1,\ldots,2g+2$, $\alpha_k(\xi)$,
$k=1,\ldots,g+1$ be a solution of the system (\ref{deformation_3}),
with the following initial conditions:
\begin{enumerate}
\item $\lambda_k(0)$ are the branch points of $Y$.
\item $\alpha_k(0)=\alpha_k$.
\end{enumerate}
(This solution is non-singular at least for sufficiently small
$\xi$). 

Let $Y(\xi)$ be a family of hyperelliptic Riemann surface with branch
points $\lambda_k(\xi)$, $\Lambda_0(\xi)=\alpha_{g+1}(\xi)$.

Then for all sufficiently small real $\xi$
\begin{enumerate}
\item The pair $(Y(\xi),\Lambda_0(\xi))$ generates periodic solutions 
  in $x$-space of the Filament equation.
\item The non-singular solutions form an open subset in the $g+1$
  dimensional variety of all solutions of the Filament Equation corresponding
  to $Y(\xi)$.
\end{enumerate}
\end{theorem}
Proof of Theorem~\ref{filaments_deformations}. By
Lemmas~\ref{rational_form} and \ref{real_flows} the Riemann surfaces
$Y(\xi)$ generate periodic solutions in $x$-space of the self-focusing NLS
with period $2\pi$. For all $\xi$ $\Lambda_0(\xi)$ is a zero of
the quasimomentum differential, which is real. By
Theorem~\ref{main_theorem} it is 
sufficient to prove that $\Lambda_0(\xi)$ is a double point of
$Y(\xi)$. Let $\mu_0(\xi)$ be one of the pre-images of
$\Lambda_0(\xi)$. $\Lambda_0(\xi)$ is a double point of
$Y(\xi)$ if and only if $p(\mu_0(\xi),\xi)=n$ or
$p(\mu_0(\xi),\xi)=n+\frac12$ where $n\in{\Bbb Z}$. 

From (\ref{deformation_1}) it follows, that 
\beq
\frac{d}{d\xi}p(\mu_0(\xi),\xi)=
\frac{\omega(\mu_0(\xi),\xi)}{d\lambda(\mu)}+
\frac{d}{d\mu}p(\mu,\xi)\left.\vphantom\int\right|_{\mu=\mu_0(\xi)}
\frac{d\mu_0(\xi)}{d\xi}\
\label{double_point_conservation}
\eeq
where the differential $\omega(\mu,\xi)$ is defined by the formula 
(\ref{omega_definition}). Since $\Lambda_0(\xi)$ coincides with one of the
zeroes of the quasimomentum differential, the second term on the
left-hand side of (\ref{double_point_conservation}) is equal to zero. Since
$c_{g+1}(\xi)\equiv0$ and thus  $\omega(\mu_0(\xi),\xi)=0$ and therefore
\beq
\frac{d}{d\xi}p(\mu_0(\xi),\xi)=0, \ p(\mu_0(\xi))=\hbox{const}. 
\eeq 
This completes the proof of the first part.
A proof of the second statement can be extracted from \cite{Sch96}.

\end{document}